# Non-Central Multivariate Chi-Square and Gamma Distributions

Thomas Royen

TH Bingen, University of Applied Sciences

e-mail: thomas.royen@t-online.de

A $(p-1)$-variate integral representation is given for the cumulative distribution function of the general $p$-variate non-central gamma distribution with a non-centrality matrix of any admissible rank. The real part of products of well known analytical functions is integrated over arguments from $(-\pi, \pi)$. To facilitate the computation, these formulas are given more detailed for $p \leq 3$. These $(p-1)$-variate integrals are also derived for the diagonal of a non-central complex Wishart matrix. Furthermore, some alternative formulas are given for the cases with an associated "one-factorial" $(p \times p)$-correlation matrix $R$, i.e. $R$ differs from a suitable diagonal matrix only by a matrix of rank 1, which holds in particular for all $(3 \times 3)$-$R$ with no vanishing correlation.

## 1. Introduction and Notation

The following notations are used: $\Sigma_{(n)}$ stands for $\Sigma_{n_1+...+n_p=n}$ with $n_1,...,n_p \in \mathbb{N}_0$. The spectral norm of any square matrix $B$ is denoted by $||B||$ and $|B|$ is its determinant. For a symmetrical matrix, $A > 0$ means positive definiteness of $A$, and $I$ or $I_p$ is always an identity matrix. The symbol $\Sigma$ also denotes a non-singular covariance matrix $(\sigma_{ik})$ with $\Sigma^{-1} = (\sigma^{ik})$, and $R = (r_{ik})$ is a non-singular correlation matrix with $R^{-1} = (r^{ik})$. The abbreviation Lt stands for "Laplace transform" or "Laplace transformation". Formulas from the NIST Handbook of Mathematical Functions are cited by HMF and their numbers.

The $p$-variate non-central chi-square distribution with $\nu$ degrees of freedom, an "associated" covariance matrix $\Sigma$ and a "non-centrality" matrix $\Delta$ (the $\chi_p^2(\nu, \Sigma, \Delta)$-distribution) is defined as the joint distribution of the diagonal elements of a non-central $W_p(\nu, \Sigma, \Delta)$-Wishart matrix. To avoid confusion with different multivariate distributions having univariate (non-central) chi-square marginal distributions, this distribution can also be called a (non-central) "Wishart chi-square distribution". Without loss of generality (w.l.o.g.) we can standardize $\Sigma$ to a correlation matrix $R$. The $(p \times \nu)$-matrix $Z$ of the corresponding Wishart matrix $ZZ^T$ has $\nu$ independent $N_p(\mu_j, R)$-normally distributed columns and $\Delta = \frac{1}{2} MM^T$ with the $(p \times \nu)$-matrix $M = (\mu_1,...,\mu_\nu)$. Therefore, rank$(\Delta) \leq \min(p, \nu)$. (In the literature, frequently $\Delta \Sigma^{-1}$ is called "non-centrality" matrix.) This distribution can be applied to power calculations for many statistical tests with several correlated chi-square test statistics, (see e.g. section 6 in [13]), but also in wireless communication at least the trivariate case occurs. For $p > 3$ the corresponding cumulative distribution function (cdf) is not easy to compute. The complexity of the $\chi_p^2(\nu, R, \Delta)$-distribution (or of

---

the slightly more general $\Gamma_p(\alpha, R, \Delta)$- distribution defined in (2.2) below) depends - apart from the dimension $p$ - essentially on two ranks, the rank $k$ of $\Delta$ and the minimal possible rank $m$ of $R - V$, where $V$ is a suitable non-singular - not necessarily positive definite (pos. def.) - diagonal matrix. For many applications, only distributions with rank$(\Delta) = 1$ are required, e.g. with identical expectation vectors $\mu_j$ in the above matrix $M$. These cases can be reduced to central $\chi^2_{p+1}(\nu, \Sigma)$ - distributions by theorem 4.1 in [13] or by a similar theorem in [8]. For $p = 2$ see formula (3.24) and (3.25) in section 3. However, a bivariate correlation matrix is always "one-factorial", i.e. $m = 1$. The $\chi^2_p(\nu, R, \Delta)$-cdf with $m = 1$ and $V > 0$ was already derived in [13] from a suitable representation of the Lt of the corresponding probability density (pdf). For the corresponding $\Gamma_p(\alpha, R, \Delta)$- cdf see [14]. These distributions are included in section 3 and extended to cases with an indefinite $V$. In particular, all $R_{3\times 3}$ with no vanishing correlation are "one-factorial" in this extended sense or limit cases with $|V| = 0$. In [2] a special $\chi^2_3(\nu, R, \Delta)$- pdf is found in formula (49) with rank$(\Delta) = 1$ and $r^{13} = 0$. This is a limit case of the distributions with a "one-factorial" $R_{3\times 3}$ with $v_{22} = 0$. A simpler form of the $\Gamma_p(\alpha, R, \Delta)$- cdf, $p \geq 3$, of such limit cases is also derived in section 3.

The most general $\Gamma_p(\alpha, R, \Delta)$- cdf is represented by a $p$ - variate integral over the cube $\mathcal{C}_p = (-\pi, \pi]^p$ in section 4 and by a $(p-1)$- variate one over $\mathcal{C}_{p-1}$ in section 5. These latter integrals are written in a more detailed form in formula (5.7) and (5.8) for $p = 2$ and $p = 3$ to facilitate the direct computation. At least for the above mentioned cases with $m = 1$ and an indefinite $V_{3\times 3}$ these bivariate integrals are supposed to be more favourable for a numerical evaluation than the corresponding formulas in section 3. A further advantage is the visibility of the integrand if $p \leq 3$. The extension with a non-central complex Wishart matrix is found in section 6.

More general integrals over $\mathcal{C}_p$ for convolutions of non-central gamma distributions were already given in [16]. In particular, for $p = 1$ an integral over $(0, \pi)$ is found there for the cdf of a pos. def. quadratic form of (non-central) Gaussian random variables. For $p = 2$ e.g. a bivariate integral is obtained for the cdf of Jensen's bivariate chi-square distribution (see also [5] and [17]). A survey on formulas for the central $\chi^2_p(\nu, R)$- distribution and some similar distributions is found in the appendix of [3].

## 2. Some Preliminaries and Special Functions

**Definition**. Any non-singular covariance matrix $\Sigma_{p\times p}$ is called "positive m-factorial" or only "m-factorial" if m is the lowest possible rank of $A$ in a representation $\Sigma = V + AA^T$ with any pos. def. diagonal matrix $V$, where $A$ may contain a mixture of real or pure imaginary columns. If m is the lowest possible rank in such a representation with a non-singular not necessarily pos. def. diagonal matrix $V$, then $\Sigma$ is called "real m-factorial".

If $\Sigma_{p\times p}$ is m-factorial with $V = W^{-2}$ then $W\Sigma W = I_p + BB^T$, where $B = WA$ has $p$ rows $b^j$ and m columns $b_\mu$ and - w.l.o.g. - $b_\mu^T b_\kappa = 0$, $\mu \neq \kappa$. Each $\Sigma_{p\times p} > 0$ is at most $(p-1)$- factorial with $W^{-2} = \lambda I_p$, where $\lambda > 0$ is the lowest eigenvalue of $\Sigma$.

The $\chi_p^2(\nu, R, \Delta)$ - pdf has the Lt

$$\hat{f}(t_1,...,t_p;\nu,R,\Delta) = |I_p + 2RT|^{-\nu/2} \operatorname{etr}(-2T(I_p + 2RT)^{-1}\Delta), \tag{2.1}$$

$T = diag(t_1,...,t_p) \geq 0, \ \nu \in \mathbb{N}, \ \Delta = \tfrac{1}{2}MM^T, \ \operatorname{rank}(\Delta) \leq \min(p,\nu).$

The $p$ - variate non-central standardized gamma - pdf (in the extended sense of Krishnamoorthy and Parthasarathy [6]) - the $\Gamma_p(\alpha, R, \Delta)$ - pdf - can be defined by the Lt

$$\hat{g}(t_1,...,t_p;\alpha,R,\Delta) = |I_p + RT|^{-\alpha} \operatorname{etr}(-T(I_p + RT)^{-1}\Delta). \tag{2.2}$$

For positive integers $\nu = 2\alpha$ this distribution differs from the $\chi_p^2(\nu, R, \Delta)$ - distribution only by a scale factor 2, but this distribution is also defined for some non-integer values $2\alpha$. The function $g(x_1,...,x_p;\alpha,R,\Delta)$ with the Lt from (2.2) is not for all $\alpha > 0$ a pdf, but additionally to $2\alpha \in \{1,...,p-2\}$, $p \geq 3$, $\operatorname{rank}(\Delta) \leq 2\alpha$, at least all real values $2\alpha \geq p-1$, $p \geq 2$, are admissible with any $\Delta \geq 0$, see e.g. [7] or [9] (and also [11] if $\operatorname{rank}(\Delta) = 1$).

Recently, in [20] the central $\Gamma_p(\alpha, R)$ - distribution was also shown to exist for all values $2\alpha > [(p-1)/2]$ without any restrictions for $R$. For an m-factorial $R$ with $WRW = I_p + BB^T$ the $\Gamma_p(\alpha, R)$ - cdf is given by

$$E\left(\prod_{j=1}^p G_\alpha(w_j^2 x_j, \tfrac{1}{2} b^j S_{2\alpha} b^{jT})\right), \ 2\alpha \in \mathbb{N} \text{ or } 2\alpha > m-1 \geq [(p-1)/2], \tag{2.3}$$

with the univariate non-central gamma cdf $G_\alpha(x, y)$ from (2.8) and the expectation referring to the $W_m(2\alpha, I_m)$ - Wishart matrix $S_{2\alpha}$, see [13] or [15]. If $B$ is real, then in (2.3) all non-integer values $2\alpha > m-1 \geq 0$ are admissible too. Hence, the $\Gamma_p(\alpha, R, \Delta)$ - distribution also exists for all real values $2\alpha \geq m$ with such m-factorial $R$ and $\operatorname{rank}(\Delta) = 1$. Formula (2.3) was also used in [19] to extend the proof of the Gaussian correlation conjecture in [18]. If $|I + RT|^{-1}$ is infinitely divisible, then $|I + RT|^{-\alpha}$ is the Lt of a $\Gamma_p(\alpha, R)$ - pdf for all $\alpha > 0$ and then formula (2.2) is the Lt of a $\Gamma_p(\alpha, R, \Delta)$ - pdf for all values $2\alpha \geq 1$ if $\operatorname{rank}(\Delta) = 1$. Sufficient and necessary conditions for $R$ entailing infinite divisibility are found in [1] and [4]. Special infinitely divisible $\Gamma_p(\alpha, R)$ - distributions arise by a one-factorial $R = W^{-2} + aa^T$ with a real column $a$.

If $R$ is real one-factorial with a real or pure imaginary column $a$ and possibly with one $w_k^2 < 0$, we have the relations

$$WRW = I_p + bb^T \text{ with } w_j^{-2} = 1 - a_j^2, \ b_j = w_j a_j, \ \beta := |WRW| = 1 + b^T b,$$

$$(WRW)^{-1} = I_p - \beta^{-1} bb^T. \tag{2.4}$$

We need the following functions: The univariate gamma density

$$g_\alpha(x) = (\Gamma(\alpha))^{-1} x^{\alpha-1} e^{-x}, \ \alpha, x > 0, \tag{2.5}$$

with shape parameter $\alpha$ and the corresponding cdf

$$G_\alpha(x) = \int_0^x g_\alpha(\xi)d\xi = P(\alpha,x) = \gamma(\alpha,x)/\Gamma(\alpha) \tag{2.6}$$

with the notations from HMF 8.2.4 for the incomplete gamma function $\gamma(\alpha,x)$, which can be extended to $\mathbb{C}$, cut along the negative real axis. Later on, we shall use the fact that $x^{-\alpha}\gamma(\alpha,x)$ is a single valued function on $\mathbb{C}$.

The non-central gamma density $g_\alpha(x,y)$ with non-centrality parameter $y$ is given by

$$g_\alpha(x,y) = e^{-y}\sum_{n=0}^\infty g_{\alpha+n}(x)\frac{y^n}{n!} = e^{-x-y}(\tfrac{x}{y})^{(\alpha-1)/2} I_{\alpha-1}(2\sqrt{xy}) = e^{-y} g_\alpha(x)\,_0F_1(\alpha;xy) \tag{2.7}$$

with the modified Bessel function $I_{\alpha-1}$ and the hypergeometric function $_0F_1$. The corresponding non-central cdf is

$$G_\alpha(x,y) = e^{-y}\sum_{n=0}^\infty G_{\alpha+n}(x)\frac{y^n}{n!}. \tag{2.8}$$

Besides,

$$G_\alpha^*(x,y) := e^y G_\alpha(x,y) = \sum_{k=1}^\infty g_{\alpha+k}(x)\,_0F_1(\alpha+k;xy) \tag{2.9}$$

will be used.

For $n \in \mathbb{N}$ we have the relations

$$G_{1/2+n}(x,y) = \tfrac{1}{2}\left(erf(\sqrt{x}+\sqrt{y}) + erf(\sqrt{x}-\sqrt{y})\right) - e^{-y}\sum_{k=1}^n g_{1/2+k}(x)\,_0F_1(\tfrac{1}{2}+k;xy), \tag{2.10}$$

$$G_{1+n}(x,y) = G_1(x,y) - e^{-y}\sum_{k=1}^n g_{1+k}(x)\,_0F_1(1+k;xy), \tag{2.11}$$

$$G_n(x,y) = 1 - Q_n(\sqrt{2y},\sqrt{2x})$$

with the Marcum $Q$-function $Q_n$ of order $n$. Furthermore, with

$$G_0(z) := \begin{cases} 0, & z<0 \\ \tfrac{1}{2}, & z=0 \\ 1, & z>0 \end{cases} \tag{2.12}$$

we have the integral representation, (see section 2 in [13]),

$$G_n(x,y) = e^{-x-y}(\tfrac{x}{y})^{n/2}\frac{1}{\pi}\int_0^\pi \frac{y\cos(n\varphi) - \sqrt{xy}\cos((n-1)\varphi)}{y - 2\sqrt{xy}\cos(\varphi) + x} e^{2\sqrt{xy}\cos(\varphi)}d\varphi + G_0(x-y),\ n\in\mathbb{N}_0,\ y>0. \tag{2.13}$$

The functions $G_\alpha(x,y)$ can be extended to $x,y \in \mathbb{C}$ too with the principal value of $x^\alpha$ if $\alpha \notin \mathbb{N}$.

Furthermore, with $x \in \mathbb{R},\ y \in \mathbb{C},$ we need the functions

$$\mathcal{G}_\alpha(x,y) = \sum_{n=0}^\infty G_{\alpha+n}(x)y^n = \begin{cases} \frac{1}{1-y}\left(G_\alpha(x) - y^{1-\alpha}e^{x(y-1)}G_\alpha(xy)\right),\ y\neq 1,\ \alpha>0 \\ \frac{1}{1-y}\left(G_{\alpha-1}(x) - y^{1-\alpha}e^{x(y-1)}G_{\alpha-1}(xy)\right),\ \alpha\geq 1, \\ xg_\alpha(x) + (1-\alpha+x)G_\alpha(x),\ y=1. \end{cases} \tag{2.14}$$

verified by Lt.

## 3. Distributions with "One-Factorial" Correlation Matrices and the Trivariate Distributions

Some earlier results are compiled here and extended by theorem 1. The $\Gamma_p(\alpha, R, \Delta)$ - cdf is given for all non-singular "real one-factorial" correlation matrices $R$ from (2.4) and any admissible rank of $\Delta \geq 0$ and additionally for the limit cases with $w_k^{-2} = 1 - a_k^2 \to 0_+$ for one index $k$. In particular, all non-singular $R_{3\times 3} = (r_{ik})$ with $r_{12} r_{13} r_{23} \neq 0$ are real one-factorial or limit cases of such matrices. We have the following subclasses: pure imaginary or real column $a$ with $\max a_j^2 < 1$, real column $a$ with one value $a_k^2 > 1$ and the limit cases with one value $a_k = 1$ (replacing $a$ by $-a$ avoids $a_k = -1$). An essential simplification is obtained for rank$(\Delta) = 1$. The cases with $\max a_j^2 < 1$ were already given in [13] and [14]. The totally different unified $(p-1)$ - variate integral representation in section 5 is suitable for feasible computations for $p \leq 3$ and might be preferred for a numerical evaluation at least for cases with $a_k^2 > 1$ or for the here not included case with one vanishing correlation.

The key for the formulas in theorem 1 for the $\Gamma_p(\alpha, R, \Delta)$ - cdf are suitable representations of the Lt from (2.2) with $R$ as in (2.4), which were already given by the formulas (2.16) and (2.17) in [13] for the corresponding $\chi_p^2(\nu, R, \Delta)$ - distribution with $\max a_j^2 < 1$. With

$$R = W^{-2} + aa^T \Leftrightarrow WRW = I_p + bb^T, \ W^2 = diag(w_1^2, ..., w_p^2), \ b_j = w_j a_j, \tag{3.1}$$

(one $w_k^2 < 0$ admissible),

$$D = R^{-1}\Delta R^{-1}, \ z_j = (1 + w_j^{-2} t_j)^{-1}, \ u_j = w_j^{-2} t_j z_j = 1 - z_j, \ Z = diag(z_1, ..., z_p), \ U = I_p - Z, \tag{3.2}$$
$$T = diag(t_1, ..., t_p) \geq 0,$$

we have the following identities

$$|I_p + RT|^{-\alpha} \text{etr}(-T(I_p + RT)^{-1}\Delta) =$$

$$|Z|^\alpha (1 + b^T Ub)^{-\alpha} \text{etr}(-W\Delta WU) \exp\left(\frac{b^T UW\Delta WUb}{1 + b^T Ub}\right) = \tag{3.3}$$

$$\exp(-a^T Da) |Z|^\alpha (1 + b^T Ub)^{-\alpha} \text{etr}(-W^{-1}DW^{-1}U) \exp\left(\frac{a^T ZDZa}{1 + b^T Ub}\right). \tag{3.4}$$

With rank$(\Delta) = 1$, $\Delta = \delta\delta^T = (\delta_i \delta_j)$, $D = dd^T = (d_i d_j)$,
these formulas become

$$|Z|^\alpha (1 + b^T Ub)^{-\alpha} \exp\left(-\sum_{j=1}^p \delta_j^2 w_j^2 u_j\right) \exp\left(\frac{(b^T UW\delta)^2}{1 + b^T Ub}\right) \tag{3.5}$$

and

$$\exp(-(a^T d)^2) |Z|^\alpha (1 + b^T Ub)^{-\alpha} \exp\left(-\sum_{j=1}^p d_j^2 w_j^{-2} u_j\right) \exp\left(\frac{(a^T Zd)^2}{1 + b^T Ub}\right), \tag{3.6}$$

respectively. With $\beta = |WRW| = 1 + b^T b$ we can also write

$$1 + b^T Ub = \beta(1 - \beta^{-1} b^T Zb). \tag{3.7}$$

If there is a negative value $w_k^{-2} = 1 - a_k^2$, we obtain a negative $\beta$ and the above formulas will be used later with $t_k > (a_k^2 - 1)^{-1}$ entailing $z_k < 0$ and $1 - \beta^{-1} b^T Z b > 0$.

The identity (3.3) follows by straightforward calculation from

$$T(I_p + RT)^{-1} = T(Z^{-1} + aa^T T)^{-1} = TZ(I_p + aa^T TZ)^{-1} = TZ\left(I_p - \frac{1}{1 + a^T TZa} aa^T TZ\right) =$$

$$TZ\left(I_p - \frac{1}{1 + b^T Ub} W^{-1} bb^T W^{-1} TZ\right).$$

For the equality of the representations in (3.3) and (3.4) only the equation

$$-a^T Da - \operatorname{tr}(W^{-1} DW^{-1} U) + \frac{a^T ZDZa}{1 + b^T Ub} = -\operatorname{tr}(W\Delta WU) + \frac{b^T UW\Delta WUb}{1 + b^T Ub} \tag{3.8}$$

has to be verified. The left hand side can be written as

$$-b^T W^{-1} DW^{-1} b - \operatorname{tr}(W^{-1} DW^{-1} U) + \frac{b^T (I-U) W^{-1} DW^{-1} (I-U) b}{1 + b^T Ub}.$$

After having inserted

$$W^{-1} DW^{-1} = (WRW)^{-1} W\Delta W (WRW)^{-1} = (I - \beta^{-1} bb^T) W\Delta W (I - \beta^{-1} bb^T)$$

the equality (3.8) can be verified by elementary calculations.

For the limit case $\varepsilon := 1 - a_k^2 \to 0_+$ the limit

$$\lim_{\varepsilon \to 0} (1 + \varepsilon t)^{-(\alpha+n)} \exp\left(-\frac{\varepsilon t}{1 + \varepsilon t} \frac{y}{\varepsilon}\right) = \exp(-ty)$$

of the Lt of $\varepsilon^{-1} g_{\alpha+n}(\varepsilon^{-1} x, \varepsilon^{-1} y)$ implies

$$\lim_{\varepsilon \to 0_+} G_{\alpha+n}(\varepsilon^{-1} x, \varepsilon^{-1} y) = \begin{cases} 0, & x < y \\ 1, & x > y \end{cases}. \tag{3.9}$$

For $a^T Da = \sum_{1 \leq i \leq j \leq p} d_{ij}^* a_i a_j$, $d_{ij}^* = \begin{cases} d_{ii}, & j = i \\ 2d_{ij}, & j > i \end{cases}$

we need

$$\frac{(a^T Da)^n}{n!} = \sum_{(2n)} d(n_1, \ldots, n_p) \prod_{j=1}^{p} a_j^{n_j}, \quad d(n_1, \ldots, n_p) = \sum \prod_{1 \leq i \leq j \leq p} \frac{1}{n_{ij}!} d_{ij}^{*n_{ij}}, \tag{3.10}$$

where the summation is taken over the $n_{ij}$ with $n_j = \sum_{i \leq j} n_{ij} + \sum_{i \geq j} n_{ji}$, $j = 1, \ldots, p$.

**Theorem 1**. The $\Gamma_p(\alpha, R, \Delta)$-cdf with the non-singular correlation matrix $R = W^{-2} + aa^T$ ($\Leftrightarrow WRW = I_p + bb^T$), $W = diag(w_1, \ldots, w_p)$ with possibly one imaginary value $w_k$, is given by

$$G(x_1, \ldots, x_p; \alpha, R, \Delta) =$$

$$\exp(-a^T Da) \int_0^\infty \sum_{n=0}^\infty \left(\sum_{(2n)} d(n_1, \ldots, n_p) \prod_{j=1}^{p} a_j^{n_j} G_{\alpha+n_j}(w_j^2 x_j, d_{jj} w_j^{-2} + b_j^2 y)\right) g_{\alpha+n}(y) dy, \tag{3.11}$$

$$\max a_j^2 < 1, \quad D = (d_{ij}) = R^{-1} \Delta R^{-1},$$

or by

$$\exp(-a^T D a) \mathrm{etr}(-DW^{-2}) \times$$

$$\int_0^\infty \sum_{n=0}^\infty \beta^{-\alpha-n} \left( \sum_{(2n)} d(n_1,...,n_p) \prod_{j=1}^p a_j^{n_j} G_{\alpha+n_j}^* (w_j^2 x_j, d_{jj} w_j^{-2} + \beta^{-1} b_j^2 y) \right) g_{\alpha+n}(y) dy, \qquad (3.12)$$

$\max a_j^2 > 1$ or $\max a_j^2 < 1$, $\beta = |WRW| = 1 + b^T b$, $G_{\alpha+n}^*(x,y) = e^y G_{\alpha+n}(x,y)$ from (2.9).

With $\mathrm{rank}(\Delta) = 1$, $\Delta = \delta\delta^T = (\delta_i \delta_j)$, $D = dd^T = (d_i d_j)$ these formulas are simplified to

$$G(x_1,...,x_p;\alpha,R,\Delta) =$$

$$\exp(-(a^T d)^2) \int_0^\infty \sum_{n=0}^\infty \frac{(2n)!}{n!} \left( \sum_{(2n)} \prod_{j=1}^p \frac{(a_j d_j)^{n_j}}{n_j!} G_{\alpha+n_j}(w_j^2 x_j, d_j^2 w_j^{-2} + b_j^2 y) \right) g_{\alpha+n}(y) dy = \qquad (3.13)$$

$$\int_0^\infty \int_0^\pi \prod_{j=1}^p G_\alpha(w_j^2 x_j, \delta_j^2 w_j^2 + b_j^2 y + 2 b_j \delta_j w_j y^{1/2} \cos(\varphi)) f_\alpha(\varphi) g_\alpha(y) d\varphi dy, \qquad (3.14)$$

$$f_\alpha(\varphi) = \left( B(\tfrac{1}{2}, \alpha - \tfrac{1}{2}) \right)^{-1} (\sin\varphi)^{2(\alpha-1)}, \quad \max a_j^2 < 1,$$

or to

$$\exp(-(a^T d)^2) \exp\left(-\sum_{j=1}^p d_j^2 w_j^{-2}\right) \times$$

$$\int_0^\infty \sum_{n=0}^\infty \beta^{-\alpha-n} \frac{(2n)!}{n!} \left( \sum_{(2n)} \prod_{j=1}^p \frac{(a_j d_j)^{n_j}}{n_j!} G_{\alpha+n_j}^*(w_j^2 x_j, d_j^2 w_j^{-2} + \beta^{-1} b_j^2 y) \right) g_{\alpha+n}(y) dy = \qquad (3.15)$$

$$\exp(-(a^T d)^2) \exp\left(-\sum_{j=1}^p d_j^2 w_j^{-2}\right) \times$$

$$\beta^{-\alpha} \int_0^\infty \int_0^\pi \prod_{j=1}^p G_\alpha^*(w_j^2 x_j, d_j^2 w_j^{-2} + \beta^{-1} b_j^2 y + 2\beta^{-1/2} b_j d_j w_j^{-1} y^{1/2} \cos(\varphi)) f_\alpha(\varphi) g_\alpha(y) d\varphi dy, \qquad (3.16)$$

$\max a_j^2 > 1$ or $\max a_j^2 < 1$.

If exactly one value $a_k^2 = 1$, then

$$G(x_1,...,x_p;\alpha,R,\Delta) =$$

$$\exp(-a^T D a) \int_0^{x_k} \sum_{n=0}^\infty \left( \sum_{(2n)} d(n_1,...,n_p) \prod_{j \neq k} a_j^{n_j} G_{\alpha+n_j}(w_j^2 x_j, d_{jj} w_j^{-2} + b_j^2 y) \right) g_{\alpha+n}(y) dy, \qquad (3.17)$$

and, if $\mathrm{rank}(\Delta) = 1$, this is simplified to

$$\exp(-(a^T d)^2) \int_0^{x_k} \sum_{n=0}^\infty \frac{(2n)!}{n!} \left( \sum_{(2n)} \prod_{j=1}^p \frac{(a_j d_j)^{n_j}}{n_j!} \cdot \prod_{j \neq k} G_{\alpha+n_j}(w_j^2 x_j, d_j^2 w_j^{-2} + b_j^2 y) \right) g_{\alpha+n}(y) dy = \qquad (3.18)$$

$$\iint_\mathcal{A} \prod_{j \neq k} G_\alpha(w_j^2 x_j, \delta_j^2 w_j^2 + b_j^2 y + 2 b_j \delta_j w_j y^{1/2} \cos(\varphi)) f_\alpha(\varphi) g_\alpha(y) dy d\varphi, \qquad (3.19)$$

$$\mathcal{A} = \left\{ (\varphi, y) \mid \varphi \in [0,\pi], y \geq 0, \delta_k^2 + y + 2\delta_k y^{1/2} \cos(\varphi) \leq x_k \right\}.$$

**Remarks.** With $y_\pm = (-\delta_k \cos(\varphi) \pm (x_k - \delta_k^2 \sin^2\varphi)^{1/2})^2$ and $\varphi_1 = \arcsin(|\delta_k|^{-1} x_k^{1/2})$ we obtain

$$\mathcal{A} = \begin{cases} \{(\varphi,y) \mid 0 \leq y \leq y_+(\varphi),\ 0 \leq \varphi \leq \pi\}, & 0 \leq \delta_k^2 < x_k \\ \{(\varphi,y) \mid y_-(\varphi) \leq y \leq y_+(\varphi),\ \pi - \varphi_1 \leq \varphi \leq \pi\}, & \delta_k^2 \geq x_k,\ \delta_k > 0 \\ \{(\varphi,y) \mid y_-(\varphi) \leq y \leq y_+(\varphi),\ 0 \leq \varphi \leq \varphi_1\}, & \delta_k^2 \geq x_k,\ \delta_k < 0 \end{cases} \qquad (3.20)$$

For $p \leq 3$ the coefficients $d(n_1,...,n_p)$ are

$$d(n_1,n_2) = \sum_{0 \leq n_{12} \leq \min(n_1,n_2)} \prod_{i \leq j} \frac{(d_{ij}^*)^{n_{ij}}}{n_{ij}!}, \quad n_i - n_{12} \text{ even}, \quad n_{ii} = \tfrac{1}{2}(n_i - n_{12}) \qquad (3.21)$$

$$d(n_1,n_2,n_3) = \sum_{n_{11}=0}^{[n_1/2]} \sum_{n_{22}=0}^{[n_2/2]} \sum_{n_{33}=0}^{[n_3/2]} \prod_{i \leq j} \frac{(d_{ij}^*)^{n_{ij}}}{n_{ij}!}, \quad n_{ij} = \tfrac{1}{2}(n_i + n_j - n_k) - (n_{ii} + n_{jj} - n_{kk}), \; i < j, \; k \neq i, j. \qquad (3.22)$$

**Proof of Theorem 1.** All the representations for the cdf in theorem 1 can be verified by the Lt of the corresponding densities. The Lt $\hat{g}_{\alpha+n}(t,y)$ of $w^2 g_{\alpha+n}(w^2 x, y)$ is

$$(1+w^{-2}t)^{-(\alpha+n)} \exp(-w^{-2}t(1+w^{-2}t)^{-1}y) = z^{\alpha+n} e^{-uy}, \; u = 1-z,$$

which also holds for $w^2 < 0$ if $t + w^2 > 0$ and consequently $z < 0$. Then, from the formulas (3.11), (3.12) we obtain the Lt in (3.3) and (3.4), from the formulas (3.13), (3.14) the Lt in (3.5), and from (3.15), (3.16) the Lt in (3.6). For the verification of the formulas (3.14), (3.16) the integral representation of the modified Bessel functions $I_{\alpha-1}$ from HMF10.32.1 is used for the integration over $\varphi$ within the corresponding integral for the Lt. For the verification of the cases with one $w_k^{-2} = 1 - a_k^2 < 0$, the variable $t_k$ of the Lt is supposed to be larger than $-w_k^2$, which is equivalent to $z_k < 0$. If for the moment the integrals for the corresponding asserted representations of the pdf are denoted by $h(x_1,...,x_p)$, we obtain with a sufficiently large $\tau > 0$ a Lt of $\exp(-\tau \Sigma x_i)h(x_1,...,x_p)$ coinciding with the Lt from (2.2) with $T + \tau I_p$, $T \geq 0$, instead of $T$. Therefore, $h(x_1,...,x_p)$ has the Lt in (2.2), which proves the corresponding representation. Finally, the formulas for the limit cases with $a_k^2 = 1$ follow from (3.9). □

The relationship between a $\chi_p^2(\nu, R, \Delta)$- pdf with $\Delta = \tfrac{1}{2} MM^T = \tfrac{1}{2} x_{p+1} \mu \mu^T$, $M_{p \times \nu} = (y_1 \mu,...,y_\nu \mu)$, $x_{p+1} = \sum_{j=1}^{\nu} y_j^2$, and a central $\chi_{p+1}^2(\nu, \Sigma)$- pdf with

$$\Sigma = \begin{pmatrix} R + \mu\mu^T & \mu \\ \mu^T & 1 \end{pmatrix}$$

is described in theorem 4.1 in [13] (for a similar theorem see [8]). For the corresponding densities, this theorem gives the relation

$$f(x_1,...,x_p;\nu,R,\Delta) = (\tfrac{1}{2} g_{\nu/2}(\tfrac{1}{2} x_{p+1}))^{-1} f(x_1,...,x_{p+1};\nu,\Sigma). \qquad (3.23)$$

Here, a slightly more general version is shown for the relation between a $\Gamma_{p+1}(\alpha, \Sigma)$- pdf and a $\Gamma_p(\alpha, \Sigma_0, \Delta)$- pdf with $\Sigma = \begin{pmatrix} \Sigma_0 + \delta\delta^T & \delta \\ \delta^T & 1 \end{pmatrix}$ and $\Delta = x_{p+1} \delta\delta^T$.

**Theorem 2.** The conditional density $(g_\alpha(x_{p+1}))^{-1} g(x_1,...,x_{p+1};\alpha,\Sigma)$ is equal to $g(x_1,...,x_p;\alpha,\Sigma_0,x_{p+1}\delta\delta^T)$.

**Proof.** With $T_p = diag(t_1,...,t_p)$ the Lt of $g_\alpha(x_{p+1})g(x_1,...,x_p;\alpha,\Sigma_0,x_{p+1}\delta\delta^T)$ is given by

$$|I_p + \Sigma_0 T_p|^{-\alpha} \int_0^\infty \exp\left(-x_{p+1}(\delta^T T_p (I_p + \Sigma_0 T_p)^{-1}\delta + t_{p+1})\right) g_\alpha(x_{p+1}) dx_{p+1} =$$

$$|I_p + \Sigma_0 T_p|^{-\alpha} \left(1 + t_{p+1} + \delta^T T_p (I_p + \Sigma_0 T_p)^{-1}\delta\right)^{-\alpha},$$

which coincides with the Lt $|I_{p+1} + \Sigma T|^{-\alpha}$ of $g(x_1,...,x_{p+1};\alpha,\Sigma)$ since

$$|I_{p+1} + \Sigma T| = (1+t_{p+1})|I_p + \Sigma_0 T_p + \delta\delta^T T_p - t_{p+1}(1+t_{p+1})^{-1}\delta\delta^T T_p| = (1+t_{p+1})|I_p + \Sigma_0 T_p + (1+t_{p+1})^{-1}\delta\delta^T T_p| =$$

$$(1+t_{p+1})|I_p + \Sigma_0 T_p||I_p + (1+t_{p+1})^{-1}(I_p + \Sigma_0 T_p)^{-1}\delta\delta^T T_p| \text{ and } |I_p + (1+t_{p+1})^{-1}(I_p + \Sigma_0 T_p)^{-1}\delta\delta^T T_p| =$$

$$1 + (1+t_{p+1})^{-1}\delta^T T_p (I_p + \Sigma_0 T_p)^{-1}\delta. \ \square$$

With the correlation $r$ in a $(2\times 2)$ - correlation matrix we also use the notation $\Gamma_2(\alpha, r, \Delta)$. With the central case from formula (3.11) and theorem 2 with

$$\Sigma_0 = \begin{pmatrix} 1 & r \\ r & 1 \end{pmatrix}$$

and a one-factorial correlation matrix

$$R = diag((1+\delta_1^2)^{-1/2}, (1+\delta_2^2)^{-1/2}, 1) \Sigma \, diag((1+\delta_1^2)^{-1/2}, (1+\delta_2^2)^{-1/2}, 1) = W^{-2} + aa^T, \ WRW = I_3 + bb^T,$$

we obtain the $\Gamma_2(\alpha, r, x_3\delta\delta^T)$ - cdf

$$G(x_1, x_2; \alpha, r, x_3\delta\delta^T) = (g(x_3))^{-1} \int_0^\infty \left(\prod_{j=1}^2 G_\alpha((1+\delta_j^2)^{-1} w_j^2 x_j, b_j^2 y)\right) w_3^2 g_\alpha(w_3^2 x_3, b_3^2 y) g_\alpha(y) dy. \tag{3.24}$$

In a similar way we can apply formula (3.12) for the central case if $\max a_j^2 > 1$.

Alternatively, we can use the absolutely convergent series for a trivariate density in formula (3.5b) in [12]. With $R^{-1} = (r^{ij})$ and $Q = (q_{ij})$, $q_{ij} = r^{ij}(r^{ii}r^{jj})^{-1/2}$, this leads to

$$G(x_1, x_2; \alpha, r, x_3\delta\delta^T) = |Q|^\alpha (g_\alpha(x_3))^{-1} \left(r^{33} g_\alpha(r^{33}x_3) \prod_{j=1}^2 G_\alpha((1+\delta_j^2)^{-1} r^{jj} x_j) + \right.$$

$$\left.\sum_{N=2}^\infty (-1)^N \sum_{m_1+m_2+m_3=M-\varepsilon_N} q(\alpha, m_1, m_2, m_3, Q) r^{33} g_{\alpha+M-m_3}(r^{33}x_3) \prod_{j=1}^2 G_{\alpha+M-m_j}((1+\delta_j^2)^{-1} r^{jj} x_j)\right), \tag{3.25}$$

$$M = [N/2], \ \varepsilon_N = N - 2M, \ q(\alpha, m_1, m_2, m_3, Q) = \frac{1}{\Gamma(\alpha)} \sum_{m=0}^{\min(m_j)} \frac{2^{2m+\varepsilon_N}\Gamma(\alpha+M-m)}{(2m+\varepsilon_N)! \prod_{j=1}^3 (m_j-m)!} \prod_{k=1}^3 q_{ij}^{2m_k+\varepsilon_N},$$

$i < j, \ k \neq i, j$.

A further alternative is the integral over $(0, \pi)$ in (5.7), which also holds for $rank(\Delta) = 2$. If theorem 2 is applied with $p = 3$, again the formulas (3.14) and (3.16) are obtained.

## 4. A p-Variate Integral Representation for the General Non-Central Distribution Function

If a $p$ - variate pdf or cdf with a known characteristic function is not explicitly available, it can be represented by a $p$ - variate integral over $\mathbb{R}^p$ by means of the Fourier inversion formula under rather general assumptions. In this section an alternative $p$ - variate integral representation for the general $\Gamma_p(\alpha,\Sigma,\Delta)$ - cdf (and the pdf) is given, which requires only integration over $(-\pi,\pi)^p$. A further advantage is the always possible reduction of this integral to a $(p-1)$ - variate one as shown in the following section. This facilitates the computation at least for $p \leq 3$ with a visible integrand. The simple idea is the representation of the required probability by an integral over a $p$ - variate Fourier series (or a $p$ - variate Laurent series).

The Lt

$$|I_p + \Sigma T|^{-\alpha} \operatorname{etr}(-T(I_p + \Sigma T)^{-1}\Delta) = \operatorname{etr}(-\Delta\Sigma^{-1})|I_p + \Sigma T|^{-\alpha} \operatorname{etr}((I_p + \Sigma T)^{-1}\Delta\Sigma^{-1}) \tag{4.1}$$

with any non-singular covariance matrix $\Sigma$ and $T = diag(t_1,...,t_p) \geq 0$ is not for all $\alpha > 0$ and all non-centrality matrices $\Delta \geq 0$ the Lt of a pdf. The following integral representation also holds for real measures whose density has the above Lt. For admissible values $\alpha$ and rank($\Delta$) of a $\Gamma_p(\alpha,\Sigma,\Delta)$ - pdf see the explanations behind (2.2) and (2.3). For a suitable form of the Lt we use the matrices

$$B = (v\Sigma)^{-1} - I_p, \quad D = v^{-1}\Sigma^{-1}\Delta\Sigma^{-1} \tag{4.2}$$

with any scale factor $v$, leading to $\|B\| < 1$, e.g.

$$v = \tfrac{1}{2}(\lambda_{min}^{-1} + \lambda_{max}^{-1}) \tag{4.3}$$

with the extremal eigenvalues of $\Sigma$, which leads to $\|B\| = (\lambda_{max} - \lambda_{min})(\lambda_{max} + \lambda_{min})^{-1}$. With

$$Z = diag(z_1,...,z_p), \quad z_j = (1+v^{-1}t_j)^{-1}, \tag{4.4}$$

we obtain

$$|I_p + \Sigma T| = |I_p + v\Sigma v^{-1}T| = |v\Sigma||(v\Sigma)^{-1} - I_p + I_p + v^{-1}T| = |v\Sigma||B + Z^{-1}| = |v\Sigma||Z^{-1}||I_p + BZ|$$

and

$$\operatorname{tr}((I_p + v\Sigma v^{-1}T)^{-1}\Delta\Sigma^{-1}) = \operatorname{tr}((B + Z^{-1})^{-1}v^{-1}\Sigma^{-1}\Delta\Sigma^{-1}) = \operatorname{tr}(Z(I_p + BZ)^{-1}D).$$

Hence, the Lt from (4.1) can be written as

$$|v\Sigma|^{-\alpha} \operatorname{etr}(-\Delta\Sigma^{-1}) |Z|^{\alpha}|I_p + BZ|^{-\alpha} \operatorname{etr}(Z(I_p + BZ)^{-1}D) = \left(\prod_{j=1}^{p} z_j^{\alpha}\right) b(z_1,...,z_p) \tag{4.5}$$

with a $p$ - variate analytical function

$$b(y_1,...,y_p) = \sum \beta(n_1,...,n_p) y_1^{n_1} \cdot ... \cdot y_p^{n_p}, \tag{4.6}$$

which is abs. convergent if $\max|y_j| < \|B\|^{-1}$.

**Theorem 3.** With the functions $\mathcal{G}_\alpha$ from (2.14) the $\Gamma_p(\alpha,\Sigma,\Delta)$ - cdf is given by

$$G(x_1,...,x_p;\alpha,\Sigma,\Delta) = |v\Sigma|^{-\alpha} \operatorname{etr}(-\Delta\Sigma^{-1}) (2\pi)^{-p} \int_{\mathcal{C}_p} \frac{\operatorname{etr}(Y^{-1}(I_p + BY^{-1})^{-1}D)}{|I_p + BY^{-1}|^{\alpha}} \prod_{j=1}^{p} \mathcal{G}_\alpha(vx_j, y_j) d\varphi_j,$$

$\mathcal{C}_p = (-\pi,\pi]^p$, $v$ from (4.3), the matrices $B$ and $D$ from (4.2), $Y = diag(y_1,...,y_p)$, $y_j = re^{i\varphi_j}$, $r > \|B\|$.

**Proof.** The asserted integral representation of the theorem is equivalent to a corresponding representation for the $\Gamma_p(\alpha, \Sigma, \Delta)$ - pdf with the derivatives $\frac{\partial}{\partial x_j} \mathcal{G}_\alpha(vx_j, y_j)$. Comparing this representation with the Lt from (4.5) we obtain

$$(2\pi)^{-p} \int_{C_p} \left( \sum \beta(n_1,...,n_p) \prod_{j=1}^{p} y_j^{-n_j} \right) \left( \sum \prod_{j=1}^{p} v g_{\alpha+n_j}(vx_j) y_j^{n_j} \right) d\varphi_1...d\varphi_p =$$

$$\sum \beta(n_1,...,n_p) \prod_{j=1}^{p} v g_{\alpha+n_j}(vx_j),$$

which is abs. convergent because of $\lim_{n\to\infty} \varepsilon^{-n} \max \left\{ \prod_{j=1}^{p} v g_{\alpha+n_j}(vx_j) \mid \sum_{j=1}^{p} n_j = n \right\} = 0$ for every $\varepsilon > 0$. The Lt of this series is $\left( \prod_{j=1}^{p} z_j^\alpha \right) \sum \beta(n_1,...,n_p) \prod_{j=1}^{p} z_j^{n_j}$, coinciding with the Lt in (4.5). Integration over $x_1,...,x_p$ concludes the proof.

**Remarks.** More generally, the scale factor $v = w^2$ can be replaced by a scale matrix $W^2 = diag(w_1^2,...,w_p^2) > 0$. Then, with

$$T_W = W^{-2}T, \; \Sigma_W = W \Sigma W, \; \Delta_W = W \Delta W, \; z_j = (1 + w_j^{-2} t_j)^{-1}, \; Z = diag(z_1,...,z_p), \; B = \Sigma_W^{-1} - I_p, \text{ and}$$
$$D = \Sigma_W^{-1} \Delta_W \Sigma_W^{-1} = W^{-1} \Sigma^{-1} \Delta \Sigma^{-1} W^{-1} \quad (4.7)$$

the Lt in (4.1) becomes

$$|W \Sigma W|^{-\alpha} \operatorname{etr}(-\Delta \Sigma^{-1}) |Z|^\alpha |I_p + BZ|^{-\alpha} \operatorname{etr}(Z(I_p + BZ)^{-1} D) = \left( \prod_{j=1}^{p} z_j^\alpha \right) \sum \beta(n_1,...,n_p) \prod_{j=1}^{p} z_j^{n_j}, \quad (4.8)$$

where the series is abs. convergent if $\max z_j < \|B\|^{-1} = \|W^{-1} \Sigma^{-1} W^{-1} - I_p\|^{-1}$, i.e. for all sufficiently large values $t_j$ if $\|B\| \geq 1$. This can occur in particular with a "natural" scaling if

$$Q = (q_{ij}) = W^{-1} \Sigma^{-1} W^{-1} \quad (4.9)$$

is a correlation matrix with $\Sigma^{-1} = (\sigma^{ij})$ and $w_j^2 = \sigma^{jj}$.

Then, with any scale matrix $W$ and the matrices $B, D$ from (4.7), the $\Gamma_p(\alpha, \Sigma, \Delta)$ - cdf is given by

$$G(x_1,...,x_p; \alpha, \Sigma, \Delta) =$$

$$|W \Sigma W|^{-\alpha} \operatorname{etr}(-\Delta \Sigma^{-1}) (2\pi)^{-p} \int_{C_p} \frac{\operatorname{etr}(Y^{-1}(I_p + BY^{-1})^{-1} D)}{|I_p + BY^{-1}|^\alpha} \prod_{j=1}^{p} \mathcal{G}_\alpha(w_j^2 x_j, y_j) d\varphi_j, \quad (4.10)$$

$$y_j = r e^{i\varphi_j}, \; r > \|B\|.$$

The proof is the same as for theorem 3, but by the series expansion and the termwise Lt the correct Lt is obtained first only for sufficiently large values $t_j$ with $\max z_j < \|B\|^{-1}$, if $\|B\| \geq 1$.

# 5. The (p-1)-Variate Integral Representation

In the integral representation from (4.10) we can perform the integration over one argument $\varphi_j$, w.l.o.g. $j = p$. With

$$B = (W\Sigma W)^{-1} - I_p = \begin{pmatrix} B_{pp} & b_p \\ b_p^T & b_{pp} \end{pmatrix}, \quad D = W^{-1}\Sigma^{-1}\Delta\Sigma^{-1}W^{-1} = \begin{pmatrix} D_{pp} & d_p \\ d_p^T & d_{pp} \end{pmatrix}, \tag{5.1}$$

$$Y_{pp} = diag(y_1, ..., y_{p-1}), \quad y_j = re^{i\varphi_j}, \quad \varphi_j \in (-\pi, \pi], \tag{5.2}$$

$$y_0 = y_0(y_1, ..., y_{p-1}) = b_p^T(Y_{pp} + B_{pp})^{-1} b_p - b_{pp}, \tag{5.3}$$

$$q = q(y_1, ..., y_{p-1}) = (b_p^T(Y_{pp} + B_{pp})^{-1}, -1) D (b_p^T(Y_{pp} + B_{pp})^{-1}, -1)^T, \tag{5.4}$$

$$L_\alpha(y_1, ..., y_{p-1}) =$$
$$etr(Y_{pp}^{-1}(I_{p-1} + B_{pp}Y_{pp}^{-1})^{-1} D_{pp}) |I_{p-1} + B_{pp}Y_{pp}^{-1}|^{-\alpha} = etr((Y_{pp} + B_{pp})^{-1} D_{pp}) |Y_{pp} + B_{pp}|^{-\alpha} |Y_{pp}|^\alpha, \tag{5.5}$$

and the functions $G_\alpha^*$ from (2.9) and $\mathcal{G}_\alpha$ from (2.14) we obtain

**Theorem 4.** With the above notations the $\Gamma_p(\alpha, \Sigma, \Delta)$ - cdf is given by

$$G(x_1, ..., x_p; \alpha, \Sigma, \Delta) =$$
$$\frac{etr(-\Delta\Sigma^{-1})}{|W\Sigma W|^\alpha} \frac{1}{(2\pi)^{p-1}} \int_{C_{p-1}} \frac{L_\alpha(y_1, ..., y_{p-1})}{(1-y_0)^\alpha} G_\alpha^*((1-y_0)w_p^2 x_p, (1-y_0)^{-1} q) \prod_{j=1}^{p-1} \mathcal{G}_\alpha(w_j^2 x_j, y_j) d\varphi_j, \quad r > \|B\|.$$

**Remark.** $\lim_{y_0 \to 1} (1-y_0)^{-\alpha} G_\alpha^*\left((1-y_0)w_p^2 x_p, \frac{q}{1-y_0}\right) = (\Gamma(\alpha+1))^{-1}(w_p^2 x_p)^\alpha {}_0F_1(\alpha+1; w_p^2 q x_p).$

For the proof we use two lemmas:

**Lemma 1.** With the notations (5.1) to (5.4) the equation

$$etr((Y+B)^{-1} D) |Y+B|^{-\alpha} = etr((Y_{pp} + B_{pp})^{-1} D_{pp}) |Y_{pp} + B_{pp}|^{-\alpha} \exp(\frac{q}{y_p - y_0})(y_p - y_0)^{-\alpha} \tag{5.6}$$

holds for $r = |y_j| > \|B\|$.

**Proof.** For

$$A := Y + B = \begin{pmatrix} A_{pp} & b_p \\ b_p^T & y_p + b_{pp} \end{pmatrix}$$

we obtain with the Schur complement

$$|A| = |A_{pp}| (y_p + b_{pp} - b_p^T A_{pp}^{-1} b_p) = |Y_{pp} + B_{pp}|(y_p - y_0),$$

$$A^{-1} = \begin{pmatrix} A_{pp}^{-1} + (y_p - y_0)^{-1} A_{pp}^{-1} b_p b_p^T A_{pp}^{-1} & -(y_p - y_0)^{-1} A_{pp}^{-1} b_p \\ -(y_p - y_0)^{-1} b_p^T A_{pp}^{-1} & (y_p - y_0)^{-1} \end{pmatrix}$$

and

$$\operatorname{tr}(A^{-1}D) = \operatorname{tr}(A_{pp}^{-1} D_{pp} + (y_p - y_0)^{-1}(A_{pp}^{-1} b_p b_p^T A_{pp}^{-1} D_{pp} - A_{pp}^{-1} b_p d_p^T)) + (y_p - y_0)^{-1}(d_{pp} - b_p^T A_{pp}^{-1} d_p) =$$

$$\operatorname{tr}(A_{pp}^{-1} D_{pp}) + (y_p - y_0)^{-1} q,$$

which implies equation (5.6). □

**Lemma 2.** Let be $q \in \mathbb{C}$ any number, $S_r = \{y \in \mathbb{C} \mid |y| = r\}$ and $y_0 \in \mathbb{C} \setminus \{1\}$ any number with $|y_0| < r = |y|$, then

$$\frac{1}{2\pi i} \oint_{S_r} \exp((y - y_0)^{-1} q) \mathcal{G}_\alpha(x, y)(y - y_0)^{-\alpha} y^{\alpha - 1} dy = (1 - y_0)^{-\alpha} G_\alpha^*((1 - y_0)x, (1 - y_0)^{-1} q).$$

**Proof.** With the binomial series we get

$$\frac{1}{2\pi i} \oint_{S_r} \frac{\partial}{\partial x} \mathcal{G}_\alpha(x, y) \frac{y^{\alpha - 1}}{(y - y_0)^{\alpha + n}} dy = \frac{1}{2\pi i} \oint_{S_r} \left( \sum_{m=0}^\infty g_{\alpha+m}(x) y^m \right) \left( \sum_{k=0}^\infty \frac{\Gamma(\alpha + n + k)}{\Gamma(\alpha + n) k!} \left( \frac{y_0}{y} \right)^k \right) y^{-n-1} dy =$$

$$\frac{1}{\Gamma(\alpha + n)} \sum_{k=0}^\infty g_{\alpha + n + k}(x) \Gamma(\alpha + n + k) y_0^k / k! = g_{\alpha + n}(x) \exp(xy_0) = (1 - y_0)^{-\alpha - n}(1 - y_0) g_{\alpha + n}((1 - y_0)x).$$

Multiplication by $q^n/n!$ and summation over n yields

$$(1 - y_0)^{-\alpha} \frac{\partial}{\partial x} G_\alpha^*((1 - y_0)x, (1 - y_0)^{-1} q)).$$

Integration over $x$ provides the assertion. □

**Proof of Theorem 4.** If $y_p$ in $Y$ is replaced by a variable $y$ with any value $|y|$, then the equation

$$|Y + B| = |Y_{pp} + B_{pp}|(y + b_{pp} - b_p^T (Y_{pp} + B_{pp})^{-1} b_p) = 0$$

always has a unique solution $y = y_0 = b_p^T (Y_{pp} + B_{pp})^{-1} b_p - b_{pp}$, where $|y_0| < r$ since $\|B_{pp}\| \leq \|B\| < r = |y_j|$, $1 \leq j \leq p - 1$. Hence, with lemma 1 and lemma 2, theorem 4 is obtained by integration over $\varphi_p$ in the integral from (4.10). □

For $p = 2$ we obtain with

$$\Sigma = \begin{pmatrix} \sigma_1^2 & \rho \sigma_1 \sigma_2 \\ \rho \sigma_1 \sigma_2 & \sigma_2^2 \end{pmatrix}, \ w_j^{-2} = \sigma_j^2 (1 - \rho^2), \ B = \begin{pmatrix} 0 & -\rho \\ -\rho & 0 \end{pmatrix}, \ D = (d_{ij}) = W^{-1} \Sigma^{-1} \Delta \Sigma^{-1} W^{-1}, \ y = re^{i\varphi}, r > |\rho|,$$

the cdf $G(x_1, x_2; \alpha, \Sigma, \Delta) = (1 - \rho^2)^\alpha \operatorname{etr}(-\Delta \Sigma^{-1}) \times$ `

$$\frac{1}{\pi} \int_0^\pi \Re \left\{ \frac{\exp(d_{11} y^{-1})}{(1 - \rho^2 y^{-1})^\alpha} G_\alpha^* \left( \frac{(1 - \rho^2 y^{-1}) x_2}{\sigma_2^2 (1 - \rho^2)}, \frac{d_{11} \rho^2 y^{-1} + 2 d_{12} \rho + d_{22} y}{y - \rho^2} \right) \mathcal{G}_\alpha \left( \frac{x_1}{\sigma_1^2 (1 - \rho^2)}, y \right) \right\} d\varphi. \quad (5.7)$$

For $p=3$ we find with $\Sigma^{-1} = (\sigma^{ij})$, $w_j^2 = \sigma^{jj}$, the correlation matrix $Q = (q_{ij}) = W^{-1}\Sigma^{-1}W^{-1}$, $B = Q - I_3$, $D = (d_{ij}) = W^{-1}\Sigma^{-1}\Delta\Sigma^{-1}W^{-1}$, $y_0(y_1, y_2) = \dfrac{q_{23}^2 y_1 + q_{13}^2 y_2 - 2q_{12}q_{13}q_{23}}{y_1 y_2 - q_{12}^2}$, $y_j = re^{i\varphi}$, $j = 1, 2$, $r > \|B\|$, $\dfrac{q(y_1, y_2)}{1 - y_0(y_1, y_2)} =$

$$\dfrac{1}{y_1 y_2 - q_{12}^2} \dfrac{(q_{13}y_2 - q_{12}q_{23}, q_{23}y_1 - q_{12}q_{13}, q_{12}^2 - y_1 y_2) D (q_{13}y_2 - q_{12}q_{23}, q_{23}y_1 - q_{12}q_{13}, q_{12}^2 - y_1 y_2)^T}{y_1 y_2 - q_{23}^2 y_1 - q_{13}^2 y_2 + 2q_{12}q_{13}q_{23} - q_{12}^2},$$

and $L_\alpha(y_1, y_2) =$

$$\exp\left(\dfrac{d_{22}y_1 + d_{11}y_2 - 2d_{12}q_{12}}{y_1 y_2 - q_{12}^2}\right)\left(1 - \dfrac{q_{12}^2}{y_1 y_2}\right)^{-\alpha}$$ the cdf

$$G(x_1, x_2, x_3; \alpha, \Sigma, \Delta) = |Q|^\alpha \, \text{etr}(-\Delta\Sigma^{-1}) \times$$

$$\dfrac{1}{2\pi^2}\int_0^\pi\int_{-\pi}^\pi \Re\left\{\dfrac{L_\alpha(y_1, y_2)}{(1 - y_0(y_1, y_2))^\alpha} G_\alpha^*\left((1 - y_0(y_1, y_2))\sigma^{33} x_3, \dfrac{q(y_1, y_2)}{1 - y_0(y_1, y_2)}\right) \prod_{j=1}^2 \mathcal{G}_\alpha\left(\sigma^{jj} x_j, y_j\right)\right\} d\varphi_1 d\varphi_2. \quad (5.8)$$

The values $(1 - \rho^2 y^{-1})^{-\alpha}$ and $(1 - q_{12}^2(y_1 y_2)^{-1})^{-\alpha}$ are defined by the binomial series and the functions $\mathcal{G}_\alpha(x, y)$ and $(1 - y_0)^{-\alpha} G_\alpha^*((1 - y_0)x, (1 - y_0)^{-1}q)$, $x > 0$, are single-valued functions. If, e.g., one element $q_{ik} = 0$, $i < k$, then $\|Q - I_3\| < 1$. If $\|Q - I_3\| > 1$, then also theorem 4 can be applied with $W^2 = vI_p$, v from (4.3), $B = v^{-1}\Sigma^{-1} - I_3$ with $\|B\| < \min(r, 1)$, and $y_0, q, L_\alpha$ from (5.3), (5.4), (5.5).

## 6. The Extension to the Distribution of the Diagonal of a Non-Central Complex Wishart Matrix

Let $X$ be a $(p \times v)$-matrix with $v$ independent $CN_p(\mu_j, \Sigma)$-distributed circular complex Gaussian column vectors $x_j$. The diagonal of $XX^*$, $(X^* = \overline{X}^T)$, has a pdf with the Lt from (4.8) with the notations from (4.7) and $v$ instead of $\alpha$, but now with a pos. def. Hermitian covariance matrix $\Sigma$ and a pos. semi-def. Hermitian "non-centrality" matrix $\Delta = MM^*$, $M = (\mu_1, ..., \mu_v)$. For the corresponding cdf $F(x_1, ..., x_p; v, \Sigma, \Delta)$ we obtain again the integral representation from (4.10), but with $\alpha$ replaced by $v$. The $(p-1)$-variate integral representation of theorem 4 can also be retained for $F(x_1, ..., x_p; v, \Sigma, \Delta)$ with the same proof if we replace within the notations (5.1), (5.3), (5.4) the row vectors $b_p^T$ and $d_p^T$ by $b_p^*$ and $d_p^*$ respectively, and $q$ by

$$q(y_1, ..., y_{p-1}) = b_p^* A_{pp}^{-1} D_{pp} A_{pp}^{-1} b_p - d_p^* A_{pp}^{-1} b_p - b_p^* A_{pp}^{-1} d_p + d_{pp}$$

with the non-Hermitian $(p-1)\times(p-1)$-matrix $A_{pp} = Y_{pp} + B_{pp}$.